\documentclass[10pt]{article}

\usepackage{amsmath,amsfonts, latexsym,graphicx}

\newcommand{\U}{{\mathcal U}}
\newcommand{\0}{{\mathbf 0}}
\newcommand{\C}{{\mathbb C}}
\newcommand{\Z}{{\mathbb Z}}

\newcommand{\Q}{{\mathbb Q}}

\newcommand{\D}{{\mathbb D}}

\newcommand{\hyp}{{\mathbb H}}

\newcommand{\arrow}[1]{\stackrel{#1}{\longrightarrow}}
\newcommand{\Adot}{\mathbf A^\bullet}

\newcommand{\Pdot}{\mathbf P^\bullet}

\newcommand{\mbf}[1]{{\mathbf #1}}

\newtheorem{defn0}{Definition}[section]
\newtheorem{prop0}[defn0]{Proposition}
\newtheorem{conj0}[defn0]{Conjecture}
\newtheorem{thm0}[defn0]{Theorem}
\newtheorem{lem0}[defn0]{Lemma}
\newtheorem{corollary0}[defn0]{Corollary}
\newtheorem{example0}[defn0]{Example}
\newtheorem{remark0}[defn0]{Remark}
\newtheorem{question0}[defn0]{Question}

\title{The Perverse Study of the Milnor Fiber}

\author{David B. Massey}

\date{}

\begin{document}

\maketitle

\begin{abstract} In this note, we provide a quick introduction to the study of the Milnor fibration via the derived category and perverse sheaves. This is primarily a dictionary for translating from the standard topological setting to the derived category and/or the Abelian category of perverse sheaves.
\end{abstract}

\sloppy

%\newpage

%\tableofcontents

%\newpage

\section{Introduction} Let $X$ be a (reduced) complex analytic space, and let $f:X\rightarrow\C$ be a complex analytic function. Suppose that $\mbf p\in f^{-1}(0)$. For convenience, we shall assume that $X$ is a (closed) complex analytic subspace of an open subset $\U$ of $\C^N$, and that $f$ is the restriction of a complex analytic function $\tilde f:\U\rightarrow C$. Locally, these assumptions are always obtainable.

By a result of L\^e, \cite{relmono}, the Milnor fibration of $f$ at $\mbf p$, introduced in \cite{milnorsing}, exists. Here, we mean the Milnor fibration inside a ball; this is to be contrasted with Milnor's initial definition of the fibration on a sphere.

Because it has no effect on the homotopy-type, the fibration is typically defined inside a {\bf closed} ball, rather than an open ball; if one then wants to examine the case inside an open ball, one simply removes the boundary of the Milnor fiber and of the total space of the Milnor fibration. To be precise, we define the {\it closed Milnor fibration} to be the following restriction of $f$:
$$
f:B_{\epsilon}(\mbf p)\cap f^{-1}(\partial \D_\delta)\rightarrow\partial\D_\delta,
$$
where $0<\delta\ll\epsilon\ll 1$, $B_\epsilon(\mbf p)$ is a closed ball of radius $\epsilon$, centered at $\mbf p$, and $\partial\D_\delta$ is the boundary circle of the disk in $\C$, of radius $\delta$, centered at $0$. The fiber of the map is the much-studied (compact) {\it Milnor fiber} $F_{f, \mbf p}$ of $f$ at $\mbf p$.

In the classical setting, where $X=\U\subseteq\C^{n+1}$ and $f=\tilde f$, the Milnor fiber $F_{f, \mbf p}$, is a compact, orientable $2n$-manifold with boundary $\partial F_{f,\mbf p}$. In fact, $F_{f, \mbf p}-\partial F_{f,\mbf p}$ is a complex $n$-dimensional manifold.

In the more-general setting, where $X$ is an arbitrary, possibly singular, analytic space, $F_{f, \mbf p}$ is a compact, $2n$-dimensional stratified space. We write $\partial F_{f,\mbf p}$ for the intersection of $F_{f,\mbf p}$ with the sphere $\partial B_{\epsilon}(\mbf p)$. Then, $F_{f, \mbf p}-\partial F_{f,\mbf p}$ is a complex $n$-dimensional stratified space.

The {\it Milnor monodromy automorphism} on the cohomology of the Milnor fiber is induced by letting the value of $f$ travel once, counterclockwise, around the circle $\partial \D_\delta$.  We denote this automorphism, in degree $i$, by
$$
T^i_{f,\mbf p}: H^i(F_{f,\mbf p}; \ \Z)\rightarrow H^i(F_{f,\mbf p}; \ \Z).
$$

\medskip

\noindent There are many, many results throughout the literature about the Milnor fiber and the Milnor monodromy; they are {\bf the} fundamental devices for studying the local topology of complex hyersurfaces.

\bigskip

Most researchers in this area are aware of the fact that the cohomology $H^*(F_{f,\mbf p}; \ \Z)$ is frequently referred to as the {\it nearby cycle} cohomology and the reduced cohomology $\widetilde H^*(F_{f,\mbf p}; \ \Z)$ is frequently referred to as the {\it vanishing cycle} cohomology. They are also aware that the nearby cycles and vanishing cycles have some description in the derived category, and that perverse sheaves arise somehow. Sadly, most papers on the subject, including our own, are long and technical, making it extraordinarily difficult to understand how the derived category statements relate to the topological situation.

This note is our attempt to provide a basic translation. We will not define the derived category, nearby or vanishing cycles, or perverse sheaves. However, hopefully, we will provide enough data to give the reader a basic dictionary, and so enable the reader to understand why the advanced machinery is very useful and powerful.  If we succeed, then this should provide motivation for further reading of more-substantive treatments of the material.

\section{Basic Notions on Complexes of Sheaves}

We continue to assume that $X$ is a complex analytic subspace of an open subset $\U$ of $\C^N$.

\smallskip

A {\it complex of sheaves} on $X$ means a chain complex of sheaves of $R$-modules on $X$, where $R$ is some reasonably nice base ring, typically $\Z$, $\Q$, or $\C$. Throughout these notes, we shall assume that our base ring is a field or a Dedekind domain (e.g., a PID).

The {\it constant complex of sheaves} $\mbf R_X^\bullet$ is usually referred to simply as the {\it constant sheaf}; it is the constant single sheaf $\mbf R_X$ in degree $0$, and zero in other degrees. We will concentrate most of our attention on the constant sheaf $\Z_X$.

Typically, we want to deal with complexes of sheaves of $R$-modules that are finite in how often their cohomological local structure changes and we also want them to be non-zero in a finite number of degrees; thus, throughout these notes, we assume, often without explicitly mentioning it, that our complexes are {\it bounded} (zero if the absolute value of the degree is big) and (cohomologically) {\it constructible} (which we will define more carefully later). Categorically, we work in what is known as the {\it derived category of bounded, constructible complexes of sheaves of $R$-modules on $X$}, where the objects are bounded, constructible complexes, but the morphisms are harder to define. However, while we shall occasionally mention that a map is a map in the derived category, one of our primary goals is to discuss, in down-to-Earth terms, important conclusions from the existence of these morphisms.

\medskip

\noindent\rule{1in}{1pt}

\medskip

Given a complex of sheaves of $\Z$-modules (or modules over another base ring) $\Adot$ on $X$, one can take the cohomology of the complex; this yields {\it cohomology sheaves}, $\mbf H^i(\Adot)$. If $\mbf p\in X$, there is a natural isomorphism between the stalk of 
$\mbf H^i(\Adot)$ at $\mbf p$ and the module obtained by first taking stalks, at $\mbf p$, in the complex $\Adot$, and then taking the cohomology of the resulting complex of modules. We denote either of these by $H^i(\Adot)_{\mbf p}$.

The integral cohomology $H^*(X; \Z)$ of $X$ is given by the {\it hypercohomology} $\hyp^*(X; \Z^\bullet_X)$. For other complexes of sheaves of $\Z$-modules $\Adot$ on $X$, the hypercohomology $\hyp^*(X; \Adot)$ is defined, and should be thought of as a generalization of $H^*(X; \Z)$ or as a generalization of sheaf cohomology. If $Y\subseteq X$, one can also consider $\hyp^*(X,Y; \Adot)$ the hypercohomology of the pair; this fits into the obvious long exact sequence with $\hyp^*(X; \Adot)$ and $\hyp^*(Y; \Adot)$.

\medskip

\noindent\rule{1in}{1pt}

\medskip

Let $j_{\mbf p}$ denote the inclusion of a point $\mbf p$ into $X$. Then, $j^*_{\mbf p}\Adot$ is the restriction of $\Adot$ to $\mbf p$; this is just the complex of modules of the stalks of $\Adot$ at $\mbf p$. Consequently, another way of writing the stalk cohomology of $\Adot$, at $\mbf p$, in degree $i$, is $H^i(j^*_{\mbf p}\Adot)$. 

Because we are assuming that $\Adot$ is a bounded, constructible complex, the stalk cohomology is also given by the hypercohomology inside a small enough open or closed ball around $\mbf p$, i.e., for $\epsilon>0$, sufficiently small, inclusions induce the isomorphisms
$$
\hyp^i(B_{\epsilon}(\mbf p)\cap X; \Adot) \ \cong \ \hyp^i(B^\circ_{\epsilon}(\mbf p)\cap X; \Adot) \ \cong \ H^i(j^*_{\mbf p}\Adot),
$$
where we have written $B^\circ_{\epsilon}(\mbf p)$ for the open ball.

\medskip

\noindent\rule{1in}{1pt}

\medskip

There is a ``dual'' operation to $j^*_{\mbf p}$, denoted $j^!_{\mbf p}$, which also produces a complex of modules (or, sheaves on $\mbf p$). The cohomology of the complex $j^!_{\mbf p}\Adot$ can be characterized as the hypercohomology of a pair: 
$$
H^i(j^!_{\mbf p}\Adot) \ \cong \ \hyp^i(B^\circ_{\epsilon}(\mbf p)\cap X, \big(B^\circ_{\epsilon}(\mbf p)-\{\mbf p\}\big)\cap X; \Adot) \ \cong \ \hyp^i(B_{\epsilon'}(\mbf p)\cap X, \partial B_{\epsilon'}(\mbf p)\cap X; \Adot),
$$
where $\epsilon>0$ is sufficiently small and $0<\epsilon'<\epsilon$.

\medskip

\noindent\rule{1in}{1pt}

\medskip

It is helpful to think about the case of the constant sheaf. We have
$$
H^i(j^*_{\mbf p}\Z^\bullet_X)\cong H^i(\mbf p; \Z)\hskip .25in\textnormal{and}\hskip .25in H^i(j^!_{\mbf p}\Z^\bullet_X)\cong H^i(B_{\epsilon}(\mbf p)\cap X, \partial B_{\epsilon}(\mbf p)\cap X; \Z).
$$
Thus, $H^i(j^*_{\mbf p}\Z^\bullet_X)=0$ if $i\neq 0$ and, for all $\mbf p\in X$, $H^0(j^*_{\mbf p}\Z^\bullet_X)\cong\Z$. 

On the other hand, using the long exact sequence of the pair $\big(B_{\epsilon}(\mbf p)\cap X, \partial B_{\epsilon}(\mbf p)\cap X\big)$,  we find that 
$$
H^i(j^!_{\mbf p}\Z^\bullet_X) \ \cong \  H^i(B_{\epsilon}(\mbf p)\cap X, \partial B_{\epsilon}(\mbf p)\cap X; \Z) \ \cong \ \widetilde H^{i-1}(\partial B_{\epsilon}(\mbf p)\cap X; \Z).
$$
Therefore, if $X$ is purely complex $n$-dimensional, then $X$ is an integral cohomology $2n$-manifold if and only if, for all $\mbf p\in X$, $H^i(j^!_{\mbf p}\Z^\bullet_X)=0$ for $i\neq 2n$, and $H^{2n}(j^!_{\mbf p}\Z^\bullet_X)\cong\Z$.

\medskip

\noindent\rule{1in}{1pt}

\medskip

There is a natural map (in the derived category) from $j^!_{\mbf p}\Adot$ to $j^*_{\mbf p}\Adot$, which, on cohomology, induces the (co)-inclusion maps from the long exact sequence of the pair:
$$
H^i(j^!_{\mbf p}\Adot) \ \cong \ \hyp^i(B^\circ_{\epsilon}(\mbf p)\cap X, \big(B^\circ_{\epsilon}(\mbf p)-\{\mbf p\}\big)\cap X; \Adot) \ \rightarrow  \  \hyp^i(B^\circ_{\epsilon}(\mbf p)\cap X; \Adot) \ \cong \ H^i(j^*_{\mbf p}\Adot).
$$

\medskip

\noindent\rule{1in}{1pt}

\medskip

For an integer $k$, the shifted complex $\Adot[k]$ is given, in degree $i$, by $\left(\Adot[k]\right)^i=\mbf A^{k+i}$, with correspondingly shifted differentials, multiplied by $(-1)^k$. Consequently, $\mbf H^i(\Adot[k])\cong \mbf H^{k+i}(\Adot)$. You should note that this ``shift by $k$'' actually may seem more like a shift by $-k$; for instance, $\Z^\bullet_X[k]$ is zero in all degrees except $-k$.

\section{The Nearby and Vanishing Cycles}

As in the introduction, let $X$ be a (closed) complex analytic subspace of an open subset $\U$ of $\C^N$, and suppose that $f:X\rightarrow\C$ is the restriction of a complex analytic function $\tilde f:\U\rightarrow C$. 

\medskip

\noindent\rule{1in}{1pt}

\medskip

Let $\Adot$ be a complex of sheaves on $X$. Then, there are two complexes of sheaves on $f^{-1}(0)$, $\psi_f\Adot$ and $\phi_f\Adot$, called the (complexes of) sheaves of {\it nearby cycles} and {\it vanishing cycles} of $\Adot$ along $f$, respectively. There are {\it monodromy automorphisms} (in the derived category) on complexes of sheaves, $T_f:\psi_f\Adot\rightarrow\psi_f\Adot$ and $\widetilde T_f:\phi_f\Adot\rightarrow\phi_f\Adot$.

\medskip

\noindent\rule{1in}{1pt}

\medskip

For $\mbf p\in f^{-1}(0)$, the stalk cohomology of the nearby cycles $\psi_f\Adot$ at $\mbf p$ gives the hypercohomology of the Milnor fiber of $f$ at $\mbf p$, with coefficients in $\Adot$. That is,
$$
H^i(\psi_f\Adot)_{\mbf p} \ \cong \ \hyp^i(F_{f,\mbf p}; \Adot).
$$
In particular, 
$$
H^i(\psi_f\Z^\bullet_X)_{\mbf p} \ \cong \ \hyp^i(F_{f,\mbf p}; \Z).
$$

The automorphism $T_f:\psi_f\Adot\rightarrow\psi_f\Adot$ yields all of the Milnor monodromy automorphisms on the stalk cohomology:
$$
T^i_{f, \mbf p}: \hyp^i(F_{f,\mbf p}; \Adot)\rightarrow  \hyp^i(F_{f,\mbf p}; \Adot).
$$

\medskip

To obtain the cohomology of the Milnor fiber modulo its boundary, one uses $j^!_{\mbf p}$; for $\mbf p\in f^{-1}(0)$,
$$
H^i(j^!_{\mbf p}\psi_f\Adot) \ \cong \ \hyp^i\left(B_{\epsilon}(\mbf p)\cap f^{-1}(0), \big(B_{\epsilon}(\mbf p)-\{\mbf p\}\big)\cap f^{-1}(0); \psi_f\Adot\right) \ \cong
$$
$$
\hyp^i\left(B_{\epsilon}(\mbf p)\cap f^{-1}(t), \partial B_{\epsilon}(\mbf p)\cap f^{-1}(t); \Adot\right) \ \cong \ \hyp^i(F_{f, \mbf p}, \partial F_{f, \mbf p}; \Adot),
$$
where $0<|t|\ll \epsilon\ll 1$.

The automorphism $T_f:\psi_f\Adot\rightarrow\psi_f\Adot$ also yields all of the relative Milnor monodromy automorphisms on the stalk cohomology:
$$
j^!_{\mbf p}T^i_{f}: \hyp^i(F_{f,\mbf p}, \partial F_{f,\mbf p}; \Adot)\rightarrow  \hyp^i(F_{f,\mbf p}, \partial F_{f,\mbf p}; \Adot).
$$

\medskip

\noindent\rule{1in}{1pt}

\medskip

For $\mbf p\in f^{-1}(0)$, the stalk cohomology of the vanishing cycles $\phi_f\Adot$ at $\mbf p$ gives the relative hypercohomology of the Milnor fiber of $f$ at $\mbf p$, with coefficients in $\Adot$. That is,
$$
H^i(\phi_f\Adot)_{\mbf p} \ \cong \ \hyp^{i+1}(B_{\epsilon}(\mbf p)\cap X, F_{f,\mbf p}; \Adot),
$$
where it is important to note the change in degrees. 

In particular, 
$$
H^i(\phi_f\Z^\bullet_X)_{\mbf p} \ \cong \ H^{i+1}(B_{\epsilon}(\mbf p)\cap X, F_{f,\mbf p}; \Z) \ \cong \ \widetilde H^i(F_{f,\mbf p}; \Z).
$$

The automorphism $\widetilde T_f:\phi_f\Adot\rightarrow\phi_f\Adot$ yields all of the ``reduced'' Milnor monodromy automorphisms on the stalk cohomology:
$$
\widetilde T^i_{f, \mbf p}: \hyp^{i+1}(B_{\epsilon}(\mbf p)\cap X, F_{f,\mbf p}; \Adot)\rightarrow \hyp^{i+1}(B_{\epsilon}(\mbf p)\cap X, F_{f,\mbf p}; \Adot).
$$

\medskip

\noindent\rule{1in}{1pt}

\medskip

There is a {\it canonical morphism} (in the derived category): 
$$
\operatorname{can}:\psi_f\Adot\rightarrow\phi_f\Adot,
$$
which, on the level of stalk cohomology, yields the coboundary map from the long exact sequence of the pair $(B_{\epsilon}(\mbf p)\cap X, F_{f,\mbf p})$:
$$
j^*_{\mbf p}({\operatorname{can}})^i:H^i(\psi_f\Adot)_{\mbf p}\cong  \hyp^i(F_{f,\mbf p}; \Adot)\rightarrow\hyp^{i+1}(B_{\epsilon}(\mbf p)\cap X, F_{f,\mbf p}; \Adot)\cong H^i(\phi_f\Adot)_{\mbf p}.
$$

On stalk cohomology, this induces the {\it canonical long exact sequence}:
$$
\cdots\rightarrow H^i(j^*_{\mbf p}\Adot)\rightarrow H^i(j^*_{\mbf p}\psi_f\Adot)\xrightarrow{ \ j^*_{\mbf p}{\operatorname{can}^i} \ }H^i(j^*_{\mbf p}\phi_f\Adot)\rightarrow H^{i+1}(j^*_{\mbf p}\Adot)\rightarrow\cdots .
$$

\medskip

\noindent\rule{1in}{1pt}

\medskip

There is a more-difficult-to-describe morphism in the other direction, from $\phi_f\Adot$ to $\psi_f\Adot$, the {\it variation morphism}:
$$
\operatorname{var}:\phi_f\Adot\rightarrow\psi_f\Adot.
$$
Rather than trying to define the variation morphism, we will, instead, give two easy-to-describe relations between the canonical and variation morphisms in the derived category :
$$
\operatorname{var}\circ\operatorname{can}=\operatorname{id}- T_f\hskip .25in\textnormal{and}\hskip .25in \operatorname{can}\circ\operatorname{var}=\operatorname{id}- \widetilde T_f.
$$
These yield the corresponding relations on stalk cohomology in each degree.

The variation induces a {\it variation long exact sequence}:
$$
\cdots\rightarrow H^{i+1}(j^!_{\mbf p}\Adot)\rightarrow H^i(j^!_{\mbf p}\phi_f\Adot)\xrightarrow{ \ j^!_{\mbf p}{\operatorname{var}^i} \ }H^i(j^!_{\mbf p}\psi_f\Adot)\rightarrow H^{i+2}(j^!_{\mbf p}\Adot)\rightarrow\cdots .
$$

\medskip

\noindent\rule{1in}{1pt}

\medskip

\section{Support, Cosupport, and Perverse Sheaves}

Suppose that we have a complex $\Adot$ of sheaves on $X$. The {\it $i$-th support of $\Adot$} is
$$
{\operatorname{supp}}^i(\Adot):= \overline{\{\mbf p\in X \ | \ H^i(j_{\mbf p}^*\Adot)\neq 0\}},
$$
where the overline denotes the topological closure. The {\it support of $\Adot$} is
$$
{\operatorname{supp}}(\Adot):= \bigcup_i{\operatorname{supp}}^i(\Adot).
$$
The {\it $i$-th cosupport of $\Adot$} is
$$
{\operatorname{cosupp}}^i(\Adot):= \overline{\{\mbf p\in X \ | \ H^i(j_{\mbf p}^!\Adot)\neq 0\}}.
$$

\medskip

\noindent\rule{1in}{1pt}

\medskip

If $\mbf p$ is an isolated point in ${\operatorname{supp}}(\Adot)$, then the natural map 
$$
j_{\mbf p}^!\Adot \ \rightarrow j_{\mbf p}^*\Adot
$$
is an isomorphism.

On stalk cohomology, one sees this as a result of the long exact hypercohomology sequence of the pair $\left(B^\circ_{\epsilon}(\mbf p)\cap X, \big(B^\circ_{\epsilon}(\mbf p)-\{\mbf p\}\big)\cap X\right)$, combined with the fact that 
$\hyp^i\left(\big(B^\circ_{\epsilon}(\mbf p)-\{\mbf p\}\big)\cap X; \Adot\right)=0$, as a result of $\mbf p$ being an isolated point in the support of $\Adot$.

\medskip

\noindent\rule{1in}{1pt}

\medskip

We wish to state an important result about the support of the vanishing cycles. We mentioned earlier that we want our complexes of sheaves $\Adot$ on $X$ to be (cohomologically) constructible. This means that there exists a Whitney stratification $\mathcal S$ of $X$ such that, for all $S\in\mathcal S$, the restriction of $\Adot$ to $S$ has locally constant cohomology and finitely-generated stalk cohomology; in this case, we say that $\Adot$ is cohomologically constructible with respect to $\mathcal S$. As a special case, note that the constant sheaf $\Z^\bullet_X$ is constructible with respect to every Whitney stratification of $X$.

Now, suppose you have a fixed Whitney stratification $\mathcal S$ of $X$, an analytic function $f:X\rightarrow\C$, and a bounded complex of sheaves $\Adot$ on $X$, which is constructible with respect to $\mathcal S$. define the {\it stratified critical locus of $f$} to be
$$
\Sigma_{\mathcal S}f \ := \ \bigcup_{S\in\mathcal S}\Sigma(f_{|_S}),
$$
where $\Sigma(f_{|_S})$ denotes the standard critical locus of an analytic function on a complex manifold. Then, 
$$
\operatorname{supp}(\phi_f\Adot) \ \subseteq \Sigma_{\mathcal S}f.
$$

In the classical situation of the constant sheaf on affine space, this is simply the well-known result that the reduced cohomology of the Milnor fiber of $f$ is zero, except possibly at critical points of $f$; this is a consequence of the Implicit Function Theorem.

\medskip

\noindent\rule{1in}{1pt}

\medskip

A constructible complex of sheaves $\Adot$ on a complex analytic space $X$ satisfies the {\it support condition} if and only if, for all integers $i$,
$$
\operatorname{dim}\big({\operatorname{supp}}^i(\Adot)\big) \ \leq \ -i,
$$
where the dimension is the complex dimension and the dimension of the empty set is, by convention, $-\infty$. Note that this means that a complex of sheaves which satisfies the support condition has all of its non-zero stalk cohomology in non-positive degrees; that is, for all $\mbf p\in X$, for all $i>0$, $H^i(j_{\mbf p}^*\Adot)=0$.

As an example, if $X$ is complex $n$-dimensional, where $n>0$, then the shifted constant sheaf $\Z^\bullet_X[n]$ satisfies the support condition, while the unshifted constant sheaf $\Z_X^\bullet$ does not.

\medskip

\noindent\rule{1in}{1pt}

\medskip

A constructible complex of sheaves $\Adot$ on a complex analytic space satisfies the {\it cosupport condition} if and only if, for all integers $i$,
$$
\operatorname{dim}\big({\operatorname{cosupp}}^i(\Adot)\big) \ \leq \ i,
$$
where, again, the dimension is the complex dimension and the dimension of the empty set is, by convention, $-\infty$. In particular, if $\Adot$ satisfies the cosupport condition, then, for all $\mbf p\in X$, for all $i<0$, $H^i(j_{\mbf p}^!\Adot)=0$.

As an example, if $X$ is a complex $n$-dimensional {\bf manifold} (or, more generally, an integral $2n$-dimensional cohomology manifold), then the shifted constant sheaf $\Z^\bullet_X[n]$ satisfies the cosupport condition. To see this, recall that we showed earlier:
$$
H^i(j^!_{\mbf p}\Z^\bullet_X) \ \cong \ \widetilde H^{i-1}(\partial B_{\epsilon}(\mbf p)\cap X; \Z).
$$
This means that
$$
H^i(j^!_{\mbf p}\Z^\bullet_X[n]) \ \cong \ H^{i+n}(j^!_{\mbf p}\Z^\bullet_X)  \ \cong \ \widetilde H^{i+n-1}(\partial B_{\epsilon}(\mbf p)\cap X; \Z).
$$
If $X$ is a complex $n$-manifold, then $\partial B_{\epsilon}(\mbf p)\cap X$ is homeomorphic to a $(2n-1)$-sphere. Thus, $H^i(j^!_{\mbf p}\Z^\bullet_X[n])=0$ if $i\neq n$, and the dimension of ${\operatorname{cosupp}}^n(\Z^\bullet_X[n])$ is precisely $n$.

\medskip

\noindent\rule{1in}{1pt}

\medskip

A {\it perverse sheaf} on a complex analytic space $X$ is a bounded, constructible complex of sheaves which satisfies both the support and cosupport conditions.

As we showed above, the shifted constant sheaf $\Z^\bullet_X[n]$ is perverse if $X$ is a complex $n$-manifold. More generally, as L\^e showed in \cite{levan}, $\Z^\bullet_X[n]$ is perverse on a purely $n$-dimensional local complete intersection.

A perverse sheaf, restricted to its support is still perverse. If $\Pdot$ is a perverse sheaf on $X$ and $\mbf p\in X$, then $H^i(\Pdot)_{\mbf p}$ is possibly non-zero only for degrees $i$ such that 
$$
-{\operatorname{dim}}_{\mbf p}\operatorname{supp}(\Pdot) \ \leq \ i \ \leq 0,
$$
where ${\operatorname{dim}}_{\mbf p}$ denotes the complex dimension at $\mbf p$. In particular, if $\mbf p$ is an isolated point in the support of a perverse sheaf $\Pdot$, then $H^i(\Pdot)_{\mbf p}$ can be non-zero only when $i=0$.

\medskip

\noindent\rule{1in}{1pt}

\medskip

Let's look at an arbitrary perverse sheaf $\Pdot$ on an analytic curve $C$ or, more precisely, on the germ of an analytic curve $C$ at a point $\mbf p$. Let $C_1,\dots, C_r$ denote the local irreducible components of $C$ at $\mbf p$. Because we care only what happens near $\mbf p$, we may assume that $\Pdot$ is constructible with respect to the Whitney stratification given by the point-stratum $\{\mbf p\}$ and $S_1:=C_1-\{\mbf p\}$, \dots, $S_r:=C_r-\{\mbf p\}$. Topologically, each of the strata $S_i$ is a punctured disk.

The perverse sheaf $\Pdot$ can have non-zero stalk cohomology only in degrees $-1$ and $0$, and 
$$
\operatorname{dim}\big({\operatorname{supp}}^0(\Adot)\big) \ \leq \ 0,
$$
i.e., $H^0$ can be non-zero at, at most, isolated points. As we are working arbitrarily close to $\mbf p$, we may thus assume that, when restricted to each $S_i$, $\Pdot$ has non-zero stalk cohomology in, at most, one degree, degree $-1$. It follows that the cohomology sheaf of the restriction $\mbf H^{-1}({\Pdot}_{|_{S_i}})$ is locally constant, and this locally constant sheaf is completely determined by the stalk at a point $\mbf q_i$ in $S_i$, and by the automorphism 
$$
h_i: H^{-1}(\Pdot)_{\mbf q_i}\rightarrow H^{-1}(\Pdot)_{\mbf q_i}
$$
induced by traveling once around the puncture.

From standard sheaf cohomology, it follows that 
$$
\hyp^{-1}(S_i; \Pdot) \ \cong \ \operatorname{ker}(\operatorname{id}-h_i)\hskip .25in \textnormal{and}\hskip .25in \hyp^{0}(S_i; \Pdot) \ \cong \ \operatorname{coker}(\textnormal{id}-h_i).
$$

\medskip

\noindent What can we say about $H^{-1}(\Pdot)_{\mbf p}$ and $H^{0}(\Pdot)_{\mbf p}$? Can they be arbitrary? No.

\medskip

Consider the hypercohomology long exact sequence of the pair $(C, C-\{\mbf p\})$:
$$
\cdots\rightarrow \hyp^{-1}(C, C-\{\mbf p\}; \Pdot)\rightarrow \hyp^{-1}(C; \Pdot)\rightarrow \hyp^{-1}(C-\{\mbf p\}; \Pdot)\rightarrow \hyp^{0}(C, C-\{\mbf p\}; \Pdot)\rightarrow 
$$
$$
\hyp^{0}(C; \Pdot)\rightarrow \hyp^{0}(C-\{\mbf p\}; \Pdot)\rightarrow \hyp^{1}(C, C-\{\mbf p\}; \Pdot)\rightarrow \hyp^{1}(C; \Pdot)\rightarrow\cdots .
$$

\medskip

\noindent Now, by the cosupport condition, $\hyp^{-1}(C, C-\{\mbf p\}; \Pdot)\cong H^{-1}(j_{\mbf p}^!\Pdot)$ equals $0$ and, by the support condition, $\hyp^{1}(C; \Pdot)\cong H^{1}(j_{\mbf p}^*\Pdot)$ equals $0$. Furthermore, from our discussion above,
$$
\hyp^{-1}(C-\{\mbf p\}; \Pdot)\cong\bigoplus_i \operatorname{ker}(\operatorname{id}-h_i)\hskip .25in \textnormal{and}\hskip .25in \hyp^{0}(C-\{\mbf p\}; \Pdot)\cong\bigoplus_i \operatorname{coker}(\operatorname{id}-h_i).
$$

We conclude that we have an exact sequence
$$
0\rightarrow H^{-1}(\Pdot)_{\mbf p}\rightarrow \bigoplus_i \operatorname{ker}(\operatorname{id}-h_i)\rightarrow H^{0}(j_{\mbf p}^!\Pdot)\rightarrow H^{0}(\Pdot)_{\mbf p}\rightarrow \bigoplus_i \operatorname{coker}(\operatorname{id}-h_i)\rightarrow H^{1}(j_{\mbf p}^!\Pdot)\rightarrow 0.
$$
Note that the inclusion of $H^{-1}(\Pdot)_{\mbf p}$ into  $\bigoplus_i \operatorname{ker}(\operatorname{id}-h_i)$ implies, in particular, that $H^{-1}(\Pdot)_{\mbf p}$ injects into $\bigoplus_i H^{-1}(\Pdot)_{\mbf q_i}$.

\section{Perverse Sheaves and Milnor Fibers}

What do perverse sheaves have to do with Milnor fibers or, in the language that we have developed, what do perverse sheaves have to with nearby and vanishing cycles?

\medskip

\noindent\rule{1in}{1pt}

\medskip

The answer is: if $\Pdot$ is a perverse sheaf on $X$, and we have a complex analytic map $f:X\rightarrow\C$, then $\psi_f\Pdot[-1]$ and $\phi_f\Pdot[-1]$ are perverse sheaves on $f^{-1}(0)$. This is frequently phrased as: ``the functors $\psi_f[-1]$ and $\phi_f[-1]$ take perverse sheaves to perverse sheaves.''

\medskip

\noindent\rule{1in}{1pt}

\medskip

Consider the classical case: $\U$ is an open subset of $\C^{n+1}$, and $f:\U\rightarrow \C$ is a complex analytic function. Then, as we have seen, $\Pdot:=\Z_{\U}^\bullet[n+1]$ is a perverse sheaf on $\U$. Furthermore, we now know that $\phi_f\Pdot[-1]$ is a perverse sheaf, whose support is contained in the critical locus $\Sigma f$ of $f$. 

Suppose that $\mbf p\in f^{-1}(0)\cap\Sigma f$, and let $s:={\textnormal{dim}}_{\mbf p}\Sigma f$. As $\phi_f\Pdot[-1]$ is perverse, we know that $H^i(\Pdot)_{\mbf p}$ is possibly non-zero only for $-s\leq i\leq 0$. What does this tell us in standard topological terms?

It tells us that
$$
H^i(\phi_f\Pdot[-1])_{\mbf p} \ = \ H^i(\phi_f\Z_{\U}^\bullet[n+1][-1])_{\mbf p} \ \cong \widetilde H^{i+n}(F_{f, \mbf p};\Z)
$$
has possibly non-zero cohomology only for $-s\leq i\leq 0$, i.e., the possibly non-trivial cohomology of the Milnor fiber occurs between degrees $n-s$ and $n$. This is the cohomological version of the classical homotopy result of Kato and Matsumoto in \cite{katomatsu}.

\medskip

\noindent\rule{1in}{1pt}

\medskip

Let's continue with $f:\U\rightarrow \C$ being a complex analytic function, but now suppose that, at a point $\mbf p\in f^{-1}(0)\cap \Sigma f$, the dimension of $\Sigma f$ is $1$. Let $C=\Sigma f$. Then, since $\operatorname{supp}\phi_f\Pdot[-1]$ is contained in $C$ and $\phi_f\Pdot[-1]$ is perverse, the restriction of $\phi_f\Pdot[-1]$ to $C$ is a perverse sheaf on a curve.

Thus, as we saw earlier, 
$$
H^{-1}(\phi_f\Pdot[-1])_p\cong \widetilde H^{n-1}(F_{f, \mbf p};\Z) \hookrightarrow \bigoplus_i \operatorname{ker}(\operatorname{id}-h_i),
$$
where $h_i:H^{n-1}(F_{f, \mbf q_i};\Z)\rightarrow H^{n-1}(F_{f, \mbf q_i};\Z)$, $\mbf q_i$ is a point near $\mbf p$ on the $i$-th irreducible germ $C_i$ of $C$ at $\mbf p$, and $h_i$ is the so-called ``vertical monodromy'' as $\mbf q_i$ moves around the punctured disk $C_i-\{\mbf p\}$. Note that, by topological triviality along Whitney strata,  $H^{n-1}(F_{f, \mbf q_i};\Z)$ is isomorphic to $\Z^{\mu_i}$, where $\mu_i$ is the Milnor number, at $\mbf q_i$, of $f$ restricted to a generic hyperplane slice through $\mbf q_i$.

\section{Verdier Dualizing}

In the derived category of bounded, constructible sheaves of $R$-modules on an complex analytic space $X$, there is a contravariant functor $\mathcal D$, the {\it Verdier dual}, such that $\mathcal D^2$ is naturally isomorphic to the identity. We do not wish to try to define $\mathcal D$, but will give some of its properties. Below, we will write $\cong$ between functors to mean ``is naturally isomorphic to''. We shall also refer to a complex analytic function $f:X\rightarrow\C$. We remind you that $\psi_f[-1]$ and $\phi_f[-1]$ are the compositions of the nearby and vanishing cycle functors with a shift by $-1$.

\begin{enumerate}

\item $\mathcal D^2\cong\textnormal{id}$;

\item $\mathcal Dj^*_{\mbf p}\cong j^!_{\mbf p}\mathcal D$ or, equivalently, $j^!_{\mbf p}\cong \mathcal Dj^*_{\mbf p}\mathcal D$ and/or $j^*_{\mbf p}\cong \mathcal Dj^!_{\mbf p}\mathcal D$;

\item $\mathcal D\big(\psi_f[-1]\big)\cong \big(\psi_f[-1]\big)\mathcal D$ and  $\mathcal D\big(\phi_f[-1]\big)\cong \big(\phi_f[-1]\big)\mathcal D$;

\item Letting $\hyp_c$ denote hypercohomology with compact supports, for every open $U \subseteq X$, there is a natural split exact sequence:
$$
0 \rightarrow \operatorname{Ext}(\hyp^{q+1}_c(U ; \Adot) , R) \rightarrow \hyp^{-q}(U ; \mathcal D\Adot) \rightarrow \operatorname{Hom}(\hyp_c^q(U ; \Adot), R) \rightarrow 0.
$$

\end{enumerate}

\medskip

\noindent\rule{1in}{1pt}

\medskip

In particular, Item (4), above tells us that, {\bf if $X$ is compact} (e.g., $X$ is a single point), then hypercohomology with compact supports is just hypercohomology, and so there is a natural split exact sequence:
$$
0 \rightarrow \operatorname{Ext}(\hyp^{q+1}(X ; \Adot) , R) \rightarrow \hyp^{-q}(X ; \mathcal D\Adot) \rightarrow \operatorname{Hom}(\hyp^q(X ; \Adot), R) \rightarrow 0.
$$
This means that we can describe the stalk cohomology of $\mathcal D\Adot$:
$$
H^{-q}(j_{\mbf p}^*\mathcal D\Adot) \ \cong \ \hyp^{-q}(\mbf p; j_{\mbf p}^*\mathcal D\Adot) \ \cong \ \hyp^{-q}(\mbf p; \mathcal Dj_{\mbf p}^!\Adot) \ \cong
$$
$$
\operatorname{Hom}\left(\hyp^q(\mbf p ; j_{\mbf p}^!\Adot), R\right)\oplus \operatorname{Ext}\left(\hyp^{q+1}(\mbf p ; j_{\mbf p}^!\Adot), R\right) \ \cong \ \operatorname{Hom}\left(H^q(j_{\mbf p}^!\Adot), R\right)\oplus \operatorname{Ext}\left(H^{q+1}(j_{\mbf p}^!\Adot), R\right).
$$

Note that, if $R$ is a field, or if $H^{q+1}(j_{\mbf p}^!\Adot)$ has no torsion, then the $\operatorname{Ext}$ term above is zero, and we obtain
$$
H^{-q}(j_{\mbf p}^*\mathcal D\Adot) \ \cong  \ \operatorname{Hom}\left(H^q(j_{\mbf p}^!\Adot), R\right).
$$
From this, it follows that, with field coefficients, the support and cosupport conditions are dual, i.e., $\Adot$ satisfies the cosupport condition if and only if  $\mathcal D\Adot$ satisfies the support condition.

\medskip

\noindent\rule{1in}{1pt}

\medskip

We say that a complex $\Adot$ is {\it self-dual} if and only if $\Adot\cong\mathcal D\Adot$. 

\medskip

If $X$ is a connected complex $n$-dimensional manifold, then $\Z^\bullet_X[n]$ is self-dual; if $X$ is also compact (and, necessarily, orientable), then the self-duality of $\Z^\bullet_X[n]$, combined with our discussions above, yields Poincar\'e duality:
$$
H^i(X; \Z) \ \cong \ \hyp^{i-n}(X; \Z^\bullet_X[n]) \ \cong \ \hyp^{i-n}\left(X; \mathcal D\big(\Z^\bullet_X[n]\big)\right) \ \cong 
$$
$$
\operatorname{Hom}\left(\hyp^{n-i}(X ;\Z^\bullet_X[n]), \Z\right)\oplus \operatorname{Ext}\left(\hyp^{n-i+1}(X ;\Z^\bullet_X[n]), \Z\right) \ \cong \ \operatorname{Hom}\left(H^{2n-i}(X ;\Z), \Z\right)\oplus \operatorname{Ext}\left(H^{2n-i+1}(X ;\Z), \Z\right).
$$

\section{The Variation Isomorphism}

The knowledgable reader may know of the {\it variation isomorphism}  $H_n(F, \partial F;\Z)\arrow{\cong} H_n(F; \Z)$ for the Milnor fiber $F$ of an isolated hypersurface singularity. Cohomologically, this becomes an isomorphism $H^n(F;\Z)\arrow{\cong} H^n(F, \partial F; \Z)$

What does this have to do with our variation morphism 
$$
\operatorname{var}:\phi_f\Adot\rightarrow\psi_f\Adot
$$
and/or the induced long exact sequences
$$
\cdots\rightarrow H^{i+1}(j^!_{\mbf p}\Adot)\rightarrow H^i(j^!_{\mbf p}\phi_f\Adot)\xrightarrow{ \ j^!_{\mbf p}{\operatorname{var}^i} \ }H^i(j^!_{\mbf p}\psi_f\Adot)\rightarrow H^{i+2}(j^!_{\mbf p}\Adot)\rightarrow\cdots ?
$$

Consider the case of a complex analytic function $f:X\rightarrow \C$ and a perverse sheaf $\Pdot$ on $X$. The ``isolated singularity'' condition is replaced with the assumption that $\mbf p$ is an isolated point in the support of $\phi_f[-1]\Pdot$.

By replacing $\Adot$ with $\Pdot$ and shifting the complexes and superscripts, the variation long exact sequence becomes
$$
\cdots\rightarrow H^{i}(j^!_{\mbf p}\Pdot)\rightarrow H^{i}(j^!_{\mbf p}\phi_f[-1]\Pdot)\xrightarrow{ \ j^!_{\mbf p}{\operatorname{var}^i} \ }H^{i}(j^!_{\mbf p}\psi_f[-1]\Pdot)\rightarrow H^{i+1}(j^!_{\mbf p}\Pdot)\rightarrow\cdots,
$$
where $\Pdot$, $\phi_f[-1]\Pdot$, and $\psi_f[-1]\Pdot$ are perverse. Recall that, by the cosupport condition $H^{i}(j^!_{\mbf p}\Pdot)=0$ if $i<0$. Also, as we are assuming that $\mbf p$ is an isolated point in the support of the perverse sheaf $\phi_f[-1]\Pdot$, $H^{i}(j^!_{\mbf p}\phi_f[-1]\Pdot)\cong H^{i}(j^*_{\mbf p}\phi_f[-1]\Pdot)$ is possibly non-zero only when $i=0$.

Therefore, the only interesting portion of our long exact sequence from above becomes
$$
0\rightarrow H^{0}(j^*_{\mbf p}\phi_f[-1]\Pdot)\xrightarrow{ \ j^!_{\mbf p}{\operatorname{var}^0} \ }H^{0}(j^!_{\mbf p}\psi_f[-1]\Pdot)\rightarrow H^{1}(j^!_{\mbf p}\Pdot)\rightarrow 0.
$$

\bigskip

Analogously, the canonical long exact sequence for $\psi_f$ and $\phi_f$ tells us that we have a short exact sequence
$$
0\rightarrow H^{-1}(j^*_{\mbf p}\Pdot)\rightarrow H^0(j^*_{\mbf p}\psi_f[-1]\Pdot)\xrightarrow{ \ j^*_{\mbf p}{\operatorname{can}^0} \ }H^0(j^*_{\mbf p}\phi_f[-1]\Pdot)\rightarrow 0.
$$

Therefore, if we assume that $H^{-1}(j^*_{\mbf p}\Pdot)=0$ and $H^{1}(j^!_{\mbf p}\Pdot)=0$, then we obtain an isomorphism $\nu_{\mbf p}$ by composing:
$$
H^0(j^*_{\mbf p}\psi_f[-1]\Pdot)\xrightarrow[\cong]{ \ j^*_{\mbf p}{\operatorname{can}^0} \ }H^0(j^*_{\mbf p}\phi_f[-1]\Pdot)\xrightarrow[\cong]{ \ j^!_{\mbf p}{\operatorname{var}^0} \ }H^{0}(j^!_{\mbf p}\psi_f[-1]\Pdot).
$$
In terms of the compact Milnor fiber $F$ of $f$ at $\mbf p$, this is an isomorphism from $\hyp^{-1}(F;\Pdot)$ to $\hyp^{-1}(F, \partial F; \Pdot)$.

In the classical case, where $\U$ is an open subset of $\C^{n+1}$, $n>0$, and  $\Pdot =\Z^\bullet_{\U}[n+1]$, it is indeed true that $H^{-1}(j^*_{\mbf p}\Pdot)=0$ and $H^{1}(j^!_{\mbf p}\Pdot)=0$, and the isomorphism $\nu_{\mbf p}$ from $H^n(F;\Z)$ to $H^n(F, \partial F;\Z)$ is the classical variation isomorphism on cohomology.

If we let $\tau$ denote the natural map $H^0(j^!_{\mbf p}\psi_f[-1]\Pdot)\rightarrow H^0(j^*_{\mbf p}\psi_f[-1]\Pdot)$, then
$$
\nu_{\mbf p}\circ\tau \ = \ \operatorname{id}-j^!_{\mbf p}T^0_{f, \mbf p}[-1],
$$
that is, the identity minus the Milnor monodromy of $\hyp^{-1}(F, \partial F; \Pdot)$. Furthermore, $$
\tau\circ\nu_{\mbf p} \ = \ \operatorname{id}-j^*_{\mbf p}T^0_{f, \mbf p}[-1],
$$
that is, the identity minus the Milnor monodromy of $\hyp^{-1}(F; \Pdot)$.
\medskip

\section{The Intersection Pairing}

Before we describe the intersection pairing in a general context, we will first describe it in the classical context.

\smallskip

Let $\U$ be an open subset of $\C^{n+1}$, where $n>0$, and let $f:\U\rightarrow \C$ be a complex analytic function. For convenience, we assume that $f(\0)=0$ and ${\textnormal{dim}}_{\0}\Sigma f=0$; thus, $f$ defines an isolated hypersurface singularity at the origin. We shall write $F$ for the compact Milnor fiber $F_{f,\0}$. Recall that $F$ is a compact, oriented, $2n$-manifold with boundary.

Then, the reduced cohomology of $F$ is concentrated in degree $n$, and 
$$
\widetilde H^n(F;\Z)\cong H^n(F;\Z)\cong  H_n(F,\partial F; \Z)\cong H^n(F,\partial F;\Z)\cong \Z^\mu,
$$
where $\mu$ is the {\it Milnor number} of $f$ at $\0$.

There is an {\it intersection pairing}
$$
<\ ,\ >: H^n(F, \partial F;\Z)\times H^n(F, \partial F;\Z)\rightarrow\Z
$$ 
given as follows:

Let $r:H^n(F, \partial F;\Z)\rightarrow H^n(F;\Z)$ be the map induced by the inclusion $(F,\emptyset)\hookrightarrow(F,\partial F)$,  let $D:H^n(F;\Z)\arrow{\cong} H_n(F, \partial F;\Z)$ be the Poincar\'e-Lefschetz duality isomorphism, and let $\tau$ be the isomorphism $H_n(F, \partial F;\Z)\arrow{\cong}\operatorname{Hom}(H^n(F, \partial F;\Z), \Z)$. Then, the intersection pairing on $H^n(F, \partial F;\Z)$ is given by
$$
<\alpha,\ \beta> \ := \ \big((\tau\circ D\circ r)(\alpha)\big)(\beta),
$$
that is, $(\tau\circ D\circ r)(\alpha)$ evaluated at $\beta$.

\medskip

\noindent\rule{1in}{1pt}

\medskip

Now, we wish to describe a generalization of this pairing in terms of stalk cohomology, and vanishing/nearby cycles. As in the previous section, we consider the case of a complex analytic function $f:X\rightarrow \C$, and a complex of sheaves $\Adot$ on $X$. We do {\bf not} assume that $\mbf p$ is an isolated point in the support of $\phi_f[-1]\Adot$.

\smallskip

\smallskip

We want to describe a pairing 
$$
H^{-k}(j^!_{\mbf p}\psi_f[-1]\Adot)\times H^{k}(j^!_{\mbf p}\psi_f[-1]\mathcal D\Adot)\rightarrow R,
$$
where $R$ is our base ring. Following our discussion in the classical case, and given all of the machinery that we have developed, this is actually quite easy.

Let 
$$
r:H^{-k}(j^!_{\mbf p}\psi_f[-1]\Adot)\rightarrow H^{-k}(j^*_{\mbf p}\psi_f[-1]\Adot)
$$ 
be the canonical map.   Let $D$ be the morphism given by the compositions

$
\displaystyle H^{-k}(j^*_{\mbf p}\psi_f[-1]\Adot)\arrow{\cong} H^{-k}(\mathcal D j^!_{\mbf p}\mathcal D\big(\psi_f[-1]\big)\Adot)\arrow{\cong} \hfill
$

$\hfill\displaystyle H^{-k}(\mathcal D j^!_{\mbf p}\big(\psi_f[-1]\big)\mathcal D\Adot)\rightarrow\operatorname{Hom}\left(H^k\left(j^!_{\mbf p}\psi_f[-1]\mathcal D\Adot\right),R\right)$

\smallskip

\noindent where, in the second isomorphism, we used that $\mathcal D$ and $\psi_f[-1]$ naturally commute, and the last map is not an isomorphism, unless $\operatorname{Ext}\left(H^{k+1}\left(j^!_{\mbf p}\psi_f[-1]\mathcal D\Adot\right), R\right)=0$.

Now, the intersection pairing on $H^{-k}(j^!_{\mbf p}\psi_f[-1]\Adot)\times H^{k}(j^!_{\mbf p}\psi_f[-1]\mathcal D\Adot)$ is given by
$$
<\alpha,\ \beta> \ := \ \big((D\circ r)(\alpha)\big)(\beta),
$$
that is, $\big((D\circ r)(\alpha)\big)$ evaluated at $\beta$.

\bigskip

This intersection pairing yields a vanishing cycle pairing: 
$$
H^{-k}(j^!_{\mbf p}\phi_f[-1]\Adot)\times H^{k}(j^!_{\mbf p}\phi_f[-1]\mathcal D\Adot)\xrightarrow{ \ (j^!_{\mbf p}{\operatorname{var}^{-k}},\  j^!_{\mbf p}{\operatorname{var}^{k}}) \ } H^{-k}(j^!_{\mbf p}\psi_f[-1]\Adot)\times H^{k}(j^!_{\mbf p}\psi_f[-1]\mathcal D\Adot)\rightarrow R.
$$

\medskip

\noindent\rule{1in}{1pt}

\medskip

Now, we will look at a more-specialized setting. Suppose that that $\Pdot$ is a perverse sheaf on $X$, that $\mbf p$ is an isolated point in the support of $\phi_f[-1]\Pdot$. Assume also that $\Pdot$ is self-dual, and fix an isomorphism $\omega:\Pdot\rightarrow \mathcal D \Pdot$. This generalizes the case of an isolated affine hypersurface singularity in an open subset $\U$ of $\C^{n+1}$, where $\Pdot$ is taken to be $\Z^\bullet_{\U}[n+1]$.

Then, the pairing that we defined above, on the vanishing cycles, is possibly non-trivial only when $k=0$, where it becomes
$$
H^{0}(j^!_{\mbf p}\phi_f[-1]\Pdot)\times H^{0}(j^!_{\mbf p}\phi_f[-1]\mathcal D\Pdot)\rightarrow R.
$$

Combined with the fact that, in our setting, $j^!_{\mbf p}\phi_f[-1]\Pdot\cong j^*_{\mbf p}\phi_f[-1]\Pdot$, and using the duality $\omega:\Pdot\rightarrow\mathcal D\Pdot$ in the second factor, we obtain the ``usual'' pairing on the vanishing cycles:
$$
H^{0}(j^*_{\mbf p}\phi_f[-1]\Pdot)\times H^{0}(j^*_{\mbf p}\phi_f[-1]\mathcal \Pdot)\arrow{\cong}H^{0}(j^!_{\mbf p}\phi_f[-1]\Pdot)\times H^{0}(j^!_{\mbf p}\phi_f[-1]\mathcal D\Pdot)\rightarrow R.
$$

\newpage
\bibliographystyle{plain}
\bibliography{Masseybib}
%\printindex
\end{document}